\documentclass[12pt,a4paper,twoside]{article}
\usepackage{amsmath,amssymb,amsfonts,amsthm}
\usepackage{times,hyperref,color}
\usepackage{graphicx}

\let\mathcal\mathscr
\newtheorem{theorem}{\bf Theorem}[section]

\newtheorem{lemma}[theorem]{\bf Lemma}
\newtheorem{corollary}[theorem]{\bf Corollary}

\usepackage{epsfig}

\begin{document}

%\large

\title{Stars versus stripes Ramsey numbers}
\author{G.R. Omidi,$^{\textrm{a},\textrm{c}}$ G. Raeisi,$^{\textrm{b},\textrm{c}}$  Z. Rahimi$^{\textrm{a}}$\vspace{1 cm}\\
{\footnotesize  $^{\textrm{a}}$Department of Mathematical Sciences, Isfahan University of Technology, Isfahan, 84156-83111, Iran}\\
{\footnotesize  $^{\textrm{b}}$Department of Mathematical Sciences, Shahrekord University, \footnotesize Shahrekord, P.O.Box 115, Iran}\\
{\footnotesize  $^{\textrm{c}}$School of Mathematics, Institute for Research in Fundamental Sciences (IPM),}\\
{\footnotesize  P.O.Box: 19395-5746, Tehran, Iran}\\
{\footnotesize{romidi@cc.iut.ac.ir, g.raeisi@math.iut.ac.ir, zahra.rahimi@math.iut.ac.ir}}}

\date {}
\maketitle

\begin{abstract}
{\footnotesize
\noindent For given simple graphs $G_1, G_2, \ldots , G_t$, the  Ramsey
number  $R(G_1, G_2, \ldots, G_t)$ is the smallest positive
integer $n$ such that if the edges of the complete graph $K_n$ are
partitioned into $t$ disjoint color classes giving $t$ graphs
$H_1,H_2,\ldots,H_t$, then at least one $H_i$  has a subgraph
isomorphic to $G_i$. In this paper, for positive integers $t_1,t_2,\ldots, t_s$ and $n_1,n_2,\ldots, n_c$ the  Ramsey number $R(S_{t_1}, S_{t_2},\ldots ,S_{t_s}, n_1K_2,n_2K_2,\ldots,n_cK_2)$  is computed, where $nK_2$ denotes a matching (stripe) of size $n$, i.e., $n$ pairwise disjoint edges and $S_{n}$ is a star with $n$ edges. This result generalizes and strengthens significantly a well-known result of Cockayne and  Lorimer and also a known result of Gy\'{a}rf\'{a}s and S\'{a}rk\"{o}zy.}
\end{abstract}

\bigskip
\section{Introduction}

In this paper, we only concerned with undirected
simple finite graphs and we follow \cite{Boundy} for terminology
and notations not defined here. For a graph $G$, we denote its
vertex set, edge set, minimum degree, maximum degree and complement graph by
$V(G)$, $E(G)$,  $\delta(G)$, $\Delta(G)$ and $\bar{G}$, respectively. If $v\in
V(G)$, we use $\deg_G{(v)}$ and $N_G(v)$ $($or simply $\deg{(v)}$ and $N(v))$ to
denote the degree and the neighbors of $v$ in $G$, respectively. Also, we use $nK_2$ to denote a matching (stripe) of size $n$, i.e., $n$ pairwise disjoint edges and as usual, a complete graph on $n$ vertices, a star with $n$ edges and a balanced complete bipartite graph on $2n$ vertices are
denoted by $K_n$, $S_{n}$ and $K_{n,n}$, respectively. In addition, for disjoint subsets $A$ and $B$ of the vertex set of a graph $G$, we use $[A,B]$ to denote the bipartite subgraph of $G$ with partite sets $A$ and $B$.

\medskip
If $G$ is a graph whose edges are colored by $c$ colors, we use $G^i$, $1\leq i\leq c$, to denote the subgraph of $G$ induced by the edges of the $i$-th color. Moreover, for a vertex $v$ of $G$, we use $\deg^i(v)$ and $N^i(v)$ to denote the degree and the neighbors of $v$ in $G^i$, respectively.

\medskip
Recall that an edge coloring of $G$ is called {\it proper} if adjacent edges are assigned different colors. The minimum number of colors for a proper edge coloring of $G$ is called  the {\it chromatic index} of $G$ and is denoted by $\chi^{\prime}(G)$. It is well known that for a bipartite graph $G$, we have $\chi^{\prime}(G)=\Delta(G)$, see \cite{Boundy}.

\medskip

Let $G,G_1,G_2,\ldots,G_c$ be given simple graphs. We write $G\rightarrow (G_1,G_2,\ldots,G_c)$, if the edges of $G$ are
partitioned into $c$ disjoint color classes giving $c$ graphs
$H_1,H_2,\ldots,H_c$, then at least one $H_i$  has a subgraph
isomorphic to $G_i$. For given simple graphs $G_1, G_2, \ldots , G_c$, the {\it multicolor Ramsey number} $R(G_1, G_2, \ldots, G_c)$ is defined as the smallest positive
integer $n$ such that $K_n\rightarrow (G_1,G_2,\ldots,G_c)$. The existence of such a positive integer is
guaranteed by the Ramsey's classical result \cite{Ramsey}. For a survey on Ramsey theory, we refer
the reader to the regularly updated survey by Radziszowski \cite{survey}.

\medskip
There is very little known about $R(G_1,G_2,\ldots,G_c)$ for $c\geq 3$,
even for very special graphs. In this paper, we consider the case that $G_i$'s are stars or stripes. The  Ramsey number of stars or stripes were investigated by several authors. The  Ramsey number of stars is determined by Burr and Roberts \cite{Ramsey number of stars} and the Ramsey number for stripes was determined by Cockayne and Lorimer \cite{stripe}. In fact they showed that $R(n_1K_2,n_2K_2,\ldots,n_cK_c)=n_1+\sum_{i=1}^{c}(n_i-1)+1$ for $n_1\geq n_2\geq \cdots \geq n_c$. In \cite{gyarfas} Gy\'{a}rf\'{a}s and S\'{a}rk\"{o}zy  determined the exact value of the Ramsey number of a star versus two stripes and then they used this result to give a positive answer to a conjecture of Schelp in an asymptotic sense. It is also worth noting
that the Ramsey number for many stars and one stripe was determined in \cite{star-stripe} as follows.

\medskip
\begin{theorem}\label{star-stripe}{\rm\cite{star-stripe}}
Let $t_1,t_2,\ldots, t_s$ be positive integers, $\Sigma=\sum_{i=1}^{s}(t_i-1)$ and $n\geq 1$.  Then\\

\noindent 1)~$R(S_{t_1},S_{t_2},\ldots, S_{t_s}, nK_2)=2n,$ if $\Sigma<n$,\\

\noindent 2)~$R(S_{t_1},S_{t_2},\ldots, S_{t_s}, nK_2)=\Sigma+n,$ if $\Sigma\geq n$, $\Sigma$ is
even and some $t_i$ is even,\\

\noindent 3)~$R(S_{t_1},S_{t_2},\ldots, S_{t_s}, nK_2)=\Sigma+n+1,$ otherwise.\\
\end{theorem}

Note that, using Theorem \ref{star-stripe} for $n=1$, we conclude that $R(S_{t_1}, S_{t_2},\ldots, S_{t_s})=\Sigma+1,$ if $\Sigma$ and at least one
$t_i$ are even and $R(S_{t_1}, S_{t_2},\ldots, S_{t_s})=\Sigma+2,$ otherwise.

\medskip
\noindent The aim of this paper is the following theorem which provides the exact value of
the  Ramsey number of any number of stars versus any number of stripes. This theorem extends known results on the Ramsey number of stars and stripes in the literature.

\medskip
\begin{theorem}\label{main3}
Let $t_1,t_2,\ldots, t_s$ and $n_1, n_2,\ldots, n_c$ be positive integers, $\Sigma=\sum_{i=1}^{s}(t_i-1)$ and $r=R(S_{t_1}, S_{t_2},\ldots ,S_{t_s})$. If $n_1\geq n_2\geq\cdots \geq n_c$, then
$$R(S_{t_1}, S_{t_2},\ldots ,S_{t_s}, n_1K_2,n_2K_2,\ldots,n_cK_2)=\max\{r+\delta,n_1\}+\sum_{i=1}^{c}(n_i-1)+1,$$
where $\delta=0$ if $\Sigma<\max\{n_1,2n_2\}$, $\Sigma$ is even and some $t_i$ is even, and $\delta=-1$, otherwise.
\end{theorem}

\medskip
As an easy corollary of Theorem \ref{main3}, we have the following result which generalizes a known result of Gy\'{a}rf\'{a}s and S\'{a}rk\"{o}zy \cite{gyarfas} on the Ramsey number of one star versus two stripes.

\medskip
\begin{corollary}\label{one star and stripes}
Let $t\geq 1$ and $n_1\geq n_2\geq\cdots \geq n_c$ be positive integers. Then $$R(S_{t}, n_1K_2,n_2K_2,\ldots,n_cK_2)=\max\{t,n_1\}
+\sum_{i=1}^{c}(n_i-1)+1.$$
\end{corollary}

\noindent By Corollary \ref{one star and stripes}, for $t\leq n_1$, $n_1=\max\{n_1,n_2,\ldots, n_c\}$, we have $$R(S_{t}, n_1K_2,n_2K_2,\ldots,n_cK_2)=R(n_1K_2,n_2K_2,\ldots,n_cK_2)=n_1+\sum_{i=1}^{c}(n_i-1)+1,$$
which strengthens significantly a well-known
result of Cockayne and Lorimer on the Ramsey number of stripes. In the other word, if $G$ is a graph obtained by deleting the edges of a graph with maximum degree $(n_1-1)$ from a complete graph on $R(n_1K_2,n_2K_2,\ldots,n_cK_2)$ vertices, then $$G\rightarrow (n_1K_2,n_2K_2,\ldots,n_cK_2).$$

In addition, we obtain the following interesting result if we investigate to Corollary \ref{one star and stripes}, when $t\geq n_1$.

\medskip
\begin{corollary}\label{cor}
Let $n_1\geq n_2\geq\cdots \geq n_c$ be arbitrary positive integers, and let $G$ be a graph on $n\geq n_1+\sum_{i=1}^{c}(n_i-1)+1$ vertices such that $\delta(G)\geq \sum_{i=1}^{c}(n_i-1)+1$. Then
$$G\rightarrow (n_1K_2,n_2K_2,\ldots,n_cK_2).$$
\end{corollary}
\medskip
\noindent\textbf{Proof. }Set $t=n-\sum_{i=1}^{c}(n_i-1)-1$. Clearly $t\geq n_1$ and so by Corollary \ref{one star and stripes}, we have $$R(S_{t}, n_1K_2,n_2K_2,\ldots,n_cK_2)=n.$$
Since $G$ has $n$ vertices and $\delta(G)\geq \sum_{i=1}^{c}(n_i-1)+1$, we have $\Delta(\bar{G})\leq t-1$, which means that $\bar{G}$ is a $S_t$-free graph and so the assertion holds by the above equation.
$\hfill\blacksquare$

\bigskip
It is also worth noting that the condition on the minimum degree in Corollary \ref{cor}
is best possible. Indeed, let $G$ be a graph on $n\geq R(n_1K_2,n_2K_2,\ldots,n_cK_2)$ vertices whose vertex set is partitioned into disjoint sets $A$, $B$ with $|A|=\Lambda=\sum_{i=1}^{c}(n_i-1)$, $|B|=n-\Lambda$ and let $E(G)=\{uv|\{u,v\}\cap A\neq\emptyset\}$. Now, set $V_0=B$ and consider a partition of vertices of $A$ into sets $V_1,V_2,\ldots ,V_c$ of sizes $n_1-1$, $n_2-1$, $\ldots,n_c-1$, respectively. For each $i=1,2,\ldots, c$ color with the $i$-th color all edges within $V_i$ or edges with one vertex in $V_i$ and one in $V_j$, where $j<i$. In this coloring, the largest monochromatic matching of color $i$ has $n_i-1$ edges, while the
minimum degree of $G$ is $\Lambda$.

\section{Proof of Theorem \ref{main3}}

In order to prove Theorem \ref{main3}, we need some lemmas. First, we start with the following simple but useful lemma.

\begin{lemma}\label{coloring bipartite graphs}
Let $t_1,t_2,\ldots,t_s$ be positive integers, $\Sigma=\sum_{i=1}^{s}(t_i-1)$ and let $H$ be a graph with $\chi^\prime(H)\leq\Sigma$. Then $E(H)$ can be decomposed into edge-disjoint subgraphs $H_1,H_2,\ldots,H_s$ such that $\Delta(H_i)\leq t_i-1.$
\end{lemma}
%%%%%%%%%%%%%%%%%%%%%%%%%%
%%%%%%%%%%%%%%%%%%%%%%%%%%
\noindent\textbf{Proof. } Consider a proper edge-coloring of $H$ with $\chi^\prime(H)$ colors. Partition the set of colors into $s$ sets $A_1,A_2,\ldots,A_s$ of sizes at most $t_1-1,t_2-1,\ldots,t_s-1$, respectively. Let $H_i$, $1\leq i\leq s$, be the subgraph of $H$ induced by the edges of colors in $A_i$. Clearly $H_i$'s are the desired subgraphs which decompose $E(H)$.
$\hfill\blacksquare$
%%%%%%%%%%%%%%%%%%%%%%%%%%%%%%%%%%%%%%%%%%%%%%%%%%%%%%%%%%%%%%%%%%%%%%%%%%%%%%%%%%%%%%%%%%%%%%%%%%%%%%
%%%%%%%%%%%%%%%%%%%%%%%%%%%%%%%%%%%%%%%%%%%%%%%%%%%%%%%%%%%%%%%%%%%%%%%%%%%%%%%%%%%%%%%%%%%%%%%%%%%%%%

\medskip
An {\it alternating cycle} in an edge colored graph is a cycle
which is properly colored i.e. no two consecutive edges in the cycle have the same color. We say
that a vertex $v$ in an edge colored graph $G$ {\it separates colors} if no component of $G-v$ is
joined to $v$ by at least two edges of different colors. Grossman and H\"{a}ggkvist gave a sufficient condition under which a two-edge colored graph must have an alternating cycle. In \cite{grossman} Grossman and H\"{a}ggkvist proved that if $G$ is a graph whose edges
are colored red and blue and there is no alternating cycle in $G$, then $G$ contains a vertex $v$ that separates the colors. Bang-Jensen and G. Gutin asked whether Grossman and H\"{a}ggkvist's result could be extended to edge-colored
graphs in general, where there is no constraint on the number of colors. In \cite{alter} Yeo gave an affirmative answer to this question as follows.

\medskip
\begin{theorem}\label{alternating}\rm(\cite{alter})
If $G$ is a $c$-edge-colored graph, $c\geq 2$, with no alternating
cycle, then there is a vertex $v\in V(G)$ such that no connected
component of $G-v$ is joined to $v$ with edges of more than one
color, i.e $G$ contains a vertex separating colors.
\end{theorem}

\medskip
Let $t_1,t_2,\ldots, t_s$ be positive integers, $\Sigma=\sum_{i=1}^{s}(t_i-1)$ and $r=R(S_{t_1}, S_{t_2},\ldots ,S_{t_s})$. Also let $n_1\geq n_2\geq\cdots \geq n_c$ be positive integers and $\Lambda=\sum_{i=1}^{c}(n_i-1)$. Set

$$
f(t_1,t_2,\ldots,t_s,n_1,n_2,\ldots,n_c)=  \left\lbrace
\begin{array}{ll}
\max\{r,n_1\}+\Lambda+1  &  \Sigma~{and~ at~ least~ one}~ t_i~{is~ even}~{and}\\
&~~~~~~~~~~~~\Sigma<\max\{n_1,2n_2\} ,\vspace{.6 cm}\\

\max\{r-1,n_1\}+\Lambda+1  & ~~~~~~~~~~  otherwise.
\end{array}
\right.
$$

\medskip
In fact, $f(t_1,t_2,\ldots,t_s,n_1,n_2,\ldots,n_c)$ is the number that we claimed is equal to the Ramsey number $R(S_{t_1},S_{t_1},\ldots, S_{t_s}, n_1K_2,n_2K_2,\ldots,n_cK_2)$ in Theorem \ref{main3}. Using these notations, we have the following lemma.

\medskip
\begin{lemma}\label{lem}
Let $t_1,t_2,\ldots, t_s$ and $n_1,n_2,\ldots,n_c$ with $n_1\geq n_2\geq\cdots \geq n_c$ be positive integers and  let $G$ be a graph on $f(t_1,t_2,\ldots,t_s,n_1,n_2,\ldots,n_c)$ vertices such that $\bar{G}\nrightarrow (S_{t_1},S_{t_1},\ldots, S_{t_s}).$ Then
$$G\rightarrow (n_1K_2,n_2K_2,\ldots,n_cK_2).$$
%there is an $s$ coloring of edges of ${\overline G}$ such that the graph induced by the $i$-th color does not contain $S_{t_i}$ as a subgraph. If edges of $G$ are arbitrary colored by colors $\beta_1, \beta_2,\ldots,\beta_c$, then for some $i$, $1\leq i\leq c$, the induced graph on edges of color $\beta_i$ contains a subgraph isomorphic to $n_iK_2$.
\end{lemma}

\medskip
\noindent{\bf Proof. }Assume that the statement of this lemma is not correct and suppose that a counterexample exists. Therefore, there are some positive integers $t_1,t_2,\ldots,t_s$ and $n_1,n_2,\ldots, n_c$ with $n_1=\max\{n_1,n_2,\ldots, n_c\}$, and a graph $G$ on  $f(t_1,t_2,\ldots,t_s,n_1,n_2,\ldots,n_c)$ vertices, such that $\bar{G}\nrightarrow (S_{t_1},S_{t_1},\ldots, S_{t_s})$ and
$G\nrightarrow (n_1K_2,n_2K_2,\ldots,n_cK_2).$ Note that $c\geq 2$, by Theorem \ref{star-stripe}.

\medskip
Among all counterexamples let $G$ be a minimal one having the maximum possible number of edges, i.e. $G$ is a graph satisfies the following conditions:\\

\noindent(a)~The number of vertices  of $G$, $f(t_1,t_2,\ldots,t_s,n_1,n_2,\ldots,n_c)$, is as small as possible. \\

\noindent (b)~Among all counterexamples satisfying (a), $G$ is a counterexample with minimum $c$, i.e. no counterexample is colored with less than $c$ colors.\\

\noindent(c)~Among all counterexamples satisfying (a) and (b), $G$ is one having the maximum possible number of edges.

\medskip
\noindent
The fact $G\nrightarrow (n_1K_2,n_2K_2,\ldots,n_cK_2)$ implies that the edges of $G$ can be colored by colors
$\beta_1, \beta_2,\ldots,\beta_c$ so that for each $i$, $1\leq i\leq c$, the induced graph on edges of color $\beta_i$ does not contain a subgraph isomorphic to $n_iK_2$. Let $G^i$ be the subgraph of $G$ induced by the edges of color $\beta_i$. As $|V(G)|\geq R(n_1k_2,n_2K_2,\ldots, n_cK_2)$, we deduce that $G$ is not a complete graph. Let $u,v$ be non-adjacent vertices in $G$. As $G$ satisfies (a), (b) and (c), $n_iK_2\subseteq G^i+uv$ (to see this, it only suffices to add the edge $uv$ to $G$ and color $uv$ by $\beta_i$ and then use the property (c) of $G$) which means that $(n_i-1)K_2\subseteq G^i$. Let $M$ be the matching of size $(n_i-1)$ in $G^i$. Since $n_iK_2\nsubseteq G^i$, we must have $N^i(u), N^i(v)\subseteq V(M)$. Moreover, the fact $n_iK_2\nsubseteq G^i$ implies that for each edge $xy\in M$, the number of edges of color $i$ between $\{x,y\}$ and $\{u,v\}$ is at most 2. Thus $\deg^i(u)+\deg^i(v)\leq 2(n_i-1)$, for each $i=1,2,\ldots,c$. Therefore, $$\deg_{G}(u)+\deg_{G}(v)=\sum_{i=1}^{c}(\deg^i(u)+\deg^i(v))\leq 2\Lambda.~~~~~~~~~~~~~~~~~~~~~~~~~~(1)$$

Since $\bar{G}\nrightarrow (S_{t_1},S_{t_1},\ldots, S_{t_s})$, there is a $s$ coloring of edges of $\bar{G}$ such that the graph induced by the $i$-th color does not contain $S_{t_i}$ as a subgraph. Thus, for every vertex $v\in V(G)$, we have $\deg_{\bar{G}}(v)\leq \Sigma$. (Indeed, if $v$ is a vertex with $\deg_{\bar{G}}(v)\geq \Sigma+1$, then the Pigeonhole principle implies that any $s$ coloring of the edges of $\bar{G}$ contains a monochromatic $S_{t_i}$ of $i$-th color with center $v$, for some $i$,  a contradiction). Therefore, $\delta(G)\geq f(t_1,t_2,\ldots,t_s,n_1,n_2,\ldots,n_c)-\Sigma-1$. An easy calculation shows that $\delta(G)\geq \Lambda+1$ unless $\Sigma\geq\max\{n_1,2n_2\}$, $\Sigma$ is even and some $t_i$ is even and in this case, we have $\delta(G)\geq \Lambda$. If $\delta(G)\geq \Lambda+1$, then for every pair of vertices $u,v$, $\deg_{G}(u)+\deg_{G}(v)\geq 2(\Lambda+1)$. Using (1), we deduce that a counterexample could not exist unless $\Sigma\geq\max\{n_1,2n_2\}$, $\Sigma$ is even and some $t_i$ is even. Therefore, hereafter we may suppose that
$\Sigma$ and at least one $t_i$ is even and $\Sigma\geq\max\{n_1,2n_2\}$. Note that, in this case we have $f(t_1,t_2,\ldots,t_s,n_1,n_2,\ldots,n_c)=\Sigma+\Lambda+1$. By (1) and the fact $\delta(G)\geq \Lambda$ we conclude that for every pair of non-adjacent vertices $u,v$ in $G$:

$$~~~~~~~~~~~~~~~~~~~~~~~~~~~~~~~~~~~~~~\deg(u)=\deg(v)=\Lambda,~~~~~~~~~~~~~~~~~~~~~~~~~~~~~~~~~~~~(2)$$
$$~~~~~~~~~~~~~~~~~~\deg^i(u)+\deg^i(v)=2(n_i-1),~i=1,2,\ldots, c~~~~~~~~~~~~~~~~~~(3)$$

\medskip

\bigskip \noindent{\bf Claim 1.} $\Sigma\leq \Lambda$.

\medskip
\noindent{\it Proof of Claim 1. }On the contrary, let $\Sigma\geq \Lambda+1$. It is easy to see that

$$R(S_{t_1}, S_{t_2},\ldots, S_{t_s}, n_1K_2,n_2K_2,\ldots,n_cK_2)\leq R(S_{t_1},S_{t_2},\ldots, S_{t_s},(\Lambda+1)K_2).$$
As $\Sigma$ and some $t_i$ are even, by Theorem \ref{star-stripe} we have
$$R(S_{t_1}, S_{t_2},\ldots, S_{t_s}, (\Lambda+1)K_2)=\Sigma+\Lambda+1,$$
which implies that $$R(S_{t_1},S_{t_2},\ldots, S_{t_s}, n_1K_2,n_2K_2,\ldots,n_cK_2)\leq \Sigma+\Lambda+1.$$
This means that $K_{N}\rightarrow (S_{t_1},S_{t_2},\ldots, S_{t_s}, n_1K_2,n_2K_2,\ldots,n_cK_2)$, where $N=\Sigma+\Lambda+1.$ Therefore $\bar{G}\rightarrow (S_{t_1},S_{t_1},\ldots, S_{t_s})$ or
$G\rightarrow (n_1K_2,n_2K_2,\ldots,n_cK_2)$, a contradiction.
$\hfill\square$

\bigskip \noindent{\bf Claim 2.} $G$ is a 2-connected graph.

\medskip
\noindent{\it Proof of Claim 2.} By Claim 1, one can easily check that $\delta(G)\geq \Lambda\geq \frac{|V(G)|}{2}$, unless $\Sigma=\Lambda$. Therefore, by the Dirac's Theorem \cite{Boundy}, $G$ is a hamiltonian graph and so a 2-connected graph unless $\Sigma=\Lambda$. Now, assume that $\Sigma=\Lambda$. In this case, $|V(G)|=2\Lambda+1$ and $\delta(G)\geq \Lambda$. Clearly, $G$ is connected (in fact the diameter of $G$ is two, since every two non-adjacent vertices have a common neighbor). If there is a cut vertex $v$ of $G$, then $G-v$ has exactly two components $G_1,G_2$ with

$$G[V(G_1)\cup\{v\}]=G[V(G_2)\cup\{v\}]= K_{\Lambda+1}.$$

\medskip
Now, we claim that all edges of $G_1\cup\{v\}$ (also $G_2\cup\{v\}$) have the same color. To see this, let $v_1$ be an arbitrary vertex of $G_1$ and let the edge $vv_1$ is colored by $\beta_{i}$, for some $i$, $1\leq i\leq c$. Since $G_1\cup\{v\}$ is a complete graph and $v_1$ is an arbitrary vertex of $G_1$, in order to show that all edges of $G_1\cup\{v\}$ have the same color, it only suffices to show that all edges of $G_1\cup\{v\}$ incident to $v_1$ are of color $\beta_{i}$. On the contrary, assume that the edge $v_1v_2$  of $G_1\cup\{v\}$ is of color $\beta_{j}$, where $j\neq i$.  Now let $M_1,M_2$ be arbitrary perfect matchings in $G_1,G_2$, respectively, where $v_1v_2\in M_1$. Therefore, $|M_1\cup M_2|=\Lambda$ and we may assume that for each $t=1,2,\ldots, c$, the matching $M_1\cup M_2$ contains exactly $n_t-1$ edges of color $\beta_{t}$, since otherwise for some $1\leq i\leq c$, $G$ has a monochromatic matching of size $n_i$ with color $\beta_{i}$, which is impossible. Set $M=(M_1\setminus \{v_1v_2\})\cup M_2\cup \{vv_1\}$. Clearly $M$ contains a monochromatic matching of size $n_i$ with color $\beta_{i}$, which is again impossible. By a similar argument, all edges of $G_2\cup\{v\}$ have the same color.
Therefore at most two colors are appeared on the edges of $G$, say $\beta_{i}$ and $\beta_{j}$ (for some $i$ and $j$). Without any loss of generality, we may assume that all edges within $G_2$ are of color $\beta_{j}$ and $j\neq 1$. As $\Lambda=\Sigma\geq \max\{n_1,2n_2\}$ and $|V(G_2)|=\Lambda$, we obtain that $|V(G_2)|\geq 2n_2\geq 2n_j$ which means that $G_2$ contains a subgraph isomorphic to $n_jK_2$ of color $\beta_{j}$, a contradiction.
$\hfill\square$

%If $c=2$, then we may suppose that $G_1$ and $G_2$ are monochromatic graphs with colors $\beta_1$ and $\beta_2$, respectively.
%Now, let  $c\geq 3$. Since $n_c\geq 2$, for two perfect matchings $M_1$ and $M_2$ of $G_1$ and $G_2$ we have $|M_1\cup M_2|=\Lambda\geq(n_i-1)+(n_j-1)+1$ and so there is a monochromatic copy of $n_iK_2$ in color $\beta_{i}$ or a monochromatic copy of $n_jK_j$ in color $\beta_{j}$, a contradiction. This contradiction shows that $G$ is a 2-connected graph.
%$\hfill\square$
\medskip
Now the analysis depends on the study of certain cycles in $G$. These are alternating cycles, colored with some colors $\beta\in\{\beta_{1},\beta_2,\ldots, \beta_{c}\}$, having no two adjacent edges of the same color. The rest of the proof is devoted to prove that an alternating cycle exists in $G$.

\bigskip \noindent{\bf Claim 3.} $G$ has an alternating cycle.

\medskip
\noindent{\it Proof of Claim 3. }On the contrary, assume that $G$ does not have an alternating cycle. Thus using Theorem \ref{alternating}, $G$ has a vertex $v$ separating colors. Since $G$ is 2-connected by Claim 2, all edges of $G$ incident to $v$ have the same color, say $\beta_{i}$. Set $G'=G\setminus \{v\}$. Note that

$$|G'|=f(t_1,t_2,\ldots,t_s,m_1,m_2,\ldots,m_c)=\Sigma+\Lambda'+1,$$
where $\Lambda'=\Lambda-1$, and  $m_1,m_2,\ldots, m_c$ are the numbers $n_1,n_2,\ldots,n_i-1,\ldots, n_c$
in the decreasing order. Clearly any $s$-coloring of the edges of $\bar{G}$ induces an $s$-coloring of the edges of $\bar{G'}$. Therefore  $\bar{G'}\nrightarrow (S_{t_1},S_{t_1},\ldots, S_{t_s}).$ From the minimality of $G$ we deduce that $G'$ has a subgraph $M$ isomorphic to $(n_i-1)K_2$ whose edges are colored by $\beta_{i}$. If the degree of $v$ as a separator vertex is at least $2n_i-1$, then there is a vertex $u\in N_G(v)$ which is unsaturated by the vertices of the matching $M$. Thus adding the edge $uv$ to the matching $M$ yields a monochromatic copy of $n_iK_2$ with color $\beta_{i}$ in $G$, which is impossible. Therefore, the proof of the claim will be completed if we prove that the degree of $v$ as a separator vertex is at least $2n_i-1$.

\medskip First let all edges of $G$ incident to $v$ have color $\beta_{i}$ and $i\geq 2$. Since $\Lambda\geq \Sigma$ by Claim 1, and also $\deg_G(v)\geq\delta(G)\geq \Lambda$, we obtain that $\deg_G(v)\geq \Sigma\geq \max\{n_1,2n_2\}\geq 2n_2\geq 2n_i$ and we are done.

\medskip
Now, let all edges of $G$ incident to $v$ have color $\beta_{1}$. Let $\deg_G(v)=2n_1-k$, for some $k\geq 2$. Note that the fact $\Sigma\geq \max\{n_1,2n_2\}$ implies that $|V(G)|=\Sigma+\Lambda+1> 2n_1$ and so $v$ is not adjacent to all vertices of $G$. Therefore, by (2) we obtain that $\deg_G(v)=2n_1-k=\Lambda$. This means that $$~~~~~~~~~~~~~~~~~~~~~~~~~~~~~~~~n_1=\sum_{i=2}^{c}(n_i-1)+k-1.~~~~~~~~~~~~~~~~~~~~~~~~~~~~~(4)$$ Since $\deg_G(v)=\Lambda$ and $|V(G)|=\Sigma+\Lambda+1$, the vertex $v$ has exactly $\Sigma$ non-neighbors in $G$. Let $S$ be the set of non-neighbors of $v$ in $G$. By (3), for every vertex $z\in S$, $\deg^1(z)=k-2$ and for $i=2,3,\ldots, c$ we have $\deg^i(z)=2n_i-2$. Since for every vertex $z\in S$, $\deg^2(z)=2n_2-2$, Equation (3) implies that the graph induced by the vertices of $S$ is a complete graph. Now, we prove that $G[S]$ contains an alternating cycle. By Theorem \ref{alternating}, $G[S]$ contains an alternating cycle unless there is a vertex which separates colors. Let $z$ be a vertex of $G[S]$ separating colors and all edges of $G[S]$ incident to $v$ have the same color, say $\beta_{i}$. If $i=1$ then $k-2\geq\deg^1_{G[S]}(z)\geq \Sigma-1$ which implies that $k-1\geq \Sigma\geq n_1$, which contradicts (4). If $i\geq 2$ then $2n_i-2\geq\deg^i_{G[S]}(z)\geq \Sigma-1$ which implies that $2n_i-1\geq \Sigma\geq 2n_2\geq 2n_i$, a contradiction. This contradiction shows that if $v$ is a vertex of $G$ separating colors and all edges of $G$ incident to $v$ have the same color $\beta_{i}$, then the degree of $v$ as a separator vertex is at least $2n_i-1$, which completes the proof of the Claim 3.
$\hfill\square$

\medskip
Now let $C$ be an alternating cycle of $G$ which has $l_i$ edges colored by $\beta_{i}$, for each $i=1,2,\ldots, c$, then it has $\sum_{i=1}^{c}l_i$ vertices. For each $i$, the $l_i$ edges of $C$ colored by $\beta_{i}$ form a subgraph isomorphic to $l_iK_2$. If $G'=V(G)\setminus V(C)$, then the number of vertices in $G'$ is $$f(t_1,t_2,\ldots,t_s,m_1,m_2,\ldots,m_c)=\Sigma+\sum_{i=1}^{c}(n_i-l_i-1)+1,$$ where $m_1,m_2,\ldots, m_c$ are the numbers $n_i-l_i$ for $i=1,\ldots,c$ in the decreasing order. As $n_1\geq m_1=\max\{m_1,m_2,\ldots,m_c\}$ and $C$ is a subgraph of $G$ which is properly colored, from the minimality of $G$ we deduce that $G'$ has a monochromatic subgraph isomorphic to $(n_i-l_i)K_2$ whose edges are colored by $\beta_{i}$, for some $1\leq i\leq c$. Combining this with a monochromatic subgraph  $l_iK_2$ of color $\beta_{i}$ in $C$, we obtain a subgraph isomorphic to $n_iK_2$ with color $\beta_{i}$ in $G$, a contradiction. This contradiction shows that this lemma is true and so the proof is completed.

$\hfill\blacksquare$
%%%%%%%%%%%%%%%%%%%%%%%%%%%%%%%%%%%%%%%%%%%%%%%%%%%%%%%%%%%%%%%%%%%%%%%%%%%%%%%%%%%%%%%%%%%%%%%%%%%%%%%%%%%
%%%%%%%%%%%%%%%%%%%%%%%%%%%%%%%%%%%%%%%%%%%%%%%%%%%%%%%%%%%%%%%%%%%%%%%%%%%%%%%%%%%%%%%%%%%%%%%%%%%%%%%%%%%

\medskip
Now, we are ready to give a proof for Theorem \ref{main3} which provides the exact value of the  Ramsey number $R(S_{t_1}, S_{t_2}, \ldots, S_{t_s}, n_1K_2, n_2K_2,\ldots,n_cK_2)$.

\bigskip
\noindent{\bf Proof of Theorem \ref{main3}}. To see that the Ramsey number can not be less than the claimed number, first consider
the case that $\Sigma<\max\{n_1,2n_2\}$, $\Sigma$ is even and some $t_i$ is even. Since $\Sigma$ and some $t_i$ are even, $r=\Sigma+1$ by Theorem \ref{star-stripe}. If $\Sigma< n_1$, then consider a partition of $n_1+\Lambda$ vertices into sets $V_1,V_2,\ldots ,V_c$ of sizes $2n_1-1$, $n_2-1$, $\ldots,n_c-1$ respectively. Color with the first color all edges which are
incident with two vertices of $V_1$ and for each $i=2, \ldots, c$ color with the $i$-th color all edges having two vertices in $V_i$ or one vertex in $V_i$ and one in $V_j$ where $j<i$. Clearly, for each $i=1, 2, \ldots, c$, the graph induced by the edges of the $i$-th color does not contain a subgraph isomorphic to $n_iK_2$.

\medskip
If $n_1\leq\Sigma<2n_2$, then partition $\Sigma+\Lambda+1$ vertices into sets $V_1=A\cup B$, $V_2=C\cup D$, $V_3, \ldots, V_c$ with $|A|=|B|=\Sigma$, $|C|=2n_2-\Sigma-1$, $|D|=n_1-n_2$ and $|V_i|=n_i-1$, $3\leq i\leq c$. Color all edges contained in $B, D$ and edges in $[B,C], [D,C], [B,D], [A,D]$ by the first color $\beta_1$, all edges contained in $A, C$ and edges in $[A,C]$ by $\beta_2$. For each $i=3, 4, \ldots, c$, color with $\beta_i$ all edges having two vertices in $V_i$ or one vertex in $V_i$ and one in $V_j$ where $j<i$. Clearly, for each $i=1, 2, \ldots, c$, the graph induced by the edges of color $\beta_i$ does not contain a subgraph isomorphic to $n_iK_2$. The remaining uncolored edges are $[A,B]$ which form a copy of $K_{\Sigma,\Sigma}$. By Lemma \ref{coloring bipartite graphs}, the edges of $K_{\Sigma,\Sigma}$ can be colored by $s$-colors $\alpha_1, \alpha_2, \ldots, \alpha_s$ such that the induced graph on edges of color $\alpha_i$, $1\leq i\leq s$, does not contain $S_{t_i}$ as a subgraph. This yields an edge coloring of the complete graph on $\max\{r,n_1\}+\Lambda$ vertices with $s+c$ colors $\alpha_1, \alpha_2, \ldots, \alpha_s$ and $\beta_1, \beta_2, \ldots, \beta_c$ such that the induced graph on edges of color $\alpha_i$, $1\leq i\leq s$, does not contain $S_{t_i}$ as a subgraph and for each $i=1, 2, \ldots, c$, the induced graph on edges of color $\beta_i$ does not contain a subgraph isomorphic to $n_iK_2$. This observation shows that if $\Sigma<\max\{n_1,2n_2\}$, $\Sigma$ is even and some $t_i$ is even, then $$R(S_{t_1}, S_{t_2}, \ldots, S_{t_s}, n_1K_2,n_2K_2,\ldots, n_cK_2)\geq\max\{r,n_1\}+\Lambda+1.$$

Now assume that the case ``$\Sigma<\max\{n_1,2n_2\}$, $\Sigma$ is even and some $t_i$ is even" does not occur. Consider a partition of $n=\max\{r-1,n_1\}+\Lambda$ vertices into sets $V_0,V_1,V_2,\ldots ,V_c$ of sizes $\max\{r-1,n_1\}, n_1-1, n_2-1, \ldots,n_c-1$ respectively. For each $i=1,2,\ldots, c$, color with $\beta_i$ all edges within $V_i$ or edges with one vertex in $V_i$ and one in $V_j$, where $j<i$. Now, if $r-1\leq n_1$, then $|V_0|=n_1$, and in this case color all edges within $V_0$ by $\beta_{1}$. In fact this is a $c$-edge coloring of $K_{n_1+\Lambda}$ that does not have a matching of size $n_i$ of color $\beta_{i}$, $1\leq i\leq c$. If $r-1>n_1$, then $|V_0|=r-1$ and so there is an edge coloring of $K_{r-1}$ with $s$ colors $\alpha_1,\ldots,\alpha_s$ without a monochromatic copy of $S_{t_i}$ of color $\alpha_i$, $1\leq i\leq s$. This yields an $(s+c)$-edge coloring of $K_n$ that does not have a monochromatic star $S_{t_i}$ with color $\alpha_i$, $1\leq i\leq s$, and no monochromatic matching of size $n_i$ in color $\beta_{i}$, $1\leq i\leq c$. Therefore $$R(S_{t_1}, S_{t_2} \ldots, S_{t_s}, n_1K_2, n_2K_2,\ldots, n_cK_2)\geq f(t_1,\ldots,t_s,n_1,\ldots,n_c).$$

\noindent To prove the other direction,  consider a complete graph on $N=f(t_1,t_2,\ldots,t_s,n_1,n_2,\ldots,n_c)$ vertices whose edges are arbitrarily colored by $s+c$ colors $\alpha_1,\alpha_2,\ldots,\alpha_s$ and $\beta_1,\beta_2,\ldots,\beta_c$. Let $G$ be the graph induced by all edges of color $\beta_1,\beta_2,\ldots,\beta_c$ in $K_N$. If for each $i$, $1\leq i\leq s$, the subgraph induced by the edges of color $\alpha_i$ in $K_N$ does not contain a copy of $S_{t_i}$, then $\bar{G}\nrightarrow (S_{t_1},S_{t_1},\ldots, S_{t_s})$ and so Lemma \ref{lem} implies that $G\rightarrow (n_1K_2,n_2K_2,\ldots,n_cK_2).$ This means that for some $i$, $1\leq i\leq c$, the subgraph of $K_N$ induced on the edges of color $\beta_i$ contains a subgraph isomorphic to $n_iK_2$, which completes the proof of the theorem. $\hfill\blacksquare$

%%%%%%%%%%%%%%%%%%%%%%%%%%%%%%%%%%%%%%%%%%%%%%%%%%%%%%%%%%%%%%%%%%%%%%%%%%%%%%%%%%%%%%%%%%%%%%%%%%%%%%
%%%%%%%%%%%%%%%%%%%%%%%%%%%%%%%%%%%%%%%%%%%%%%%%%%%%%%%%%%%%%%%%%%%%%%%%%%%%

\medskip
\subsection*{Acknowledgment}
The research of the first and second authors are partially
carried out in the IPM-Isfahan Branch and in part supported respectively
by grants No. 94050217 and No. 94050057, from IPM.

\end{document}